\documentclass[a4paper,12pt,twoside]{article}
\usepackage{amsmath,amssymb,ifthen}
\usepackage[amsmath,thmmarks]{ntheorem}
\usepackage{paper}
\usepackage[all]{xy}
\usepackage{mathrsfs}
\SelectTips{cm}{11}
\DeclareMathOperator{\Spf}{Spf}
\DeclareMathOperator{\Nrd}{Nrd}
\DeclareMathOperator{\Fix}{Fix}
\DeclareMathOperator{\vol}{vol}
\DeclareMathOperator{\PGL}{PGL}
\DeclareMathOperator{\GSp}{GSp}
\DeclareMathOperator{\Lie}{Lie}
\DeclareMathOperator{\Ad}{Ad}
\DeclareMathOperator{\JL}{\mathit{JL}}
\DeclareMathOperator{\ch}{ch}

\allowdisplaybreaks

\title{Lefschetz trace formula and $\ell$-adic cohomology of Lubin-Tate tower}
\author{Yoichi Mieda}

\begin{document}

\maketitle

\begin{firstfootnote}
 Faculty of Mathematics, Kyushu University, 744 Motooka, Nishi-ku, Fukuoka-city, Fukuoka, 819--0395 Japan

 E-mail address: \texttt{mieda@math.kyushu-u.ac.jp}

 2010 \textit{Mathematics Subject Classification}.
 Primary: 22E50;
 Secondary: 14G35, 11F70.
\end{firstfootnote}

\begin{abstract}
 In this article, we investigate the alternating sum of the $\ell$-adic cohomology of the Lubin-Tate tower by the Lefschetz trace formula.
 Our method gives slightly stronger results than in the preceding work of Strauch.
\end{abstract}

\section{Introduction}
Let $F$ be a $p$-adic field, i.e., a finite extension of $\Q_p$,
and $\mathcal{O}$ the ring of integers of $F$.
For an integer $d\ge 1$, we consider the universal deformation space of formal $\mathcal{O}$-modules of height $d$.
It is called the Lubin-Tate space. By adding Drinfeld level structures,
we get a tower of affine formal schemes over the Lubin-Tate space, which is called the Lubin-Tate tower.

Let $i$ be a non-negative integer. Using the $i$th compactly supported cohomology
of the rigid generic fiber of the Lubin-Tate tower,
we obtain a representation $H^i_{\mathrm{LT}}$ of $\GL_d(F)\times D^\times\times W_F$,
where $D$ denotes the central division algebra over $F$ with invariant $1/d$,
and $W_F$ denotes the Weil group of $F$ (for the construction of $H^i_{\mathrm{LT}}$,
see Section \ref{sec:Lubin-Tate}).
The representation $H^i_{\mathrm{LT}}$ is very interesting;
in fact, it is known that for a supercuspidal representation $\pi$ of $\GL_d(F)$ the $\pi$-isotypic component
of the virtual representation $H_{\mathrm{LT}}=\sum_i(-1)^iH^i_{\mathrm{LT}}$ can be described by
the local Langlands correspondence and the local Jacquet-Langlands correspondence.
This was known as non-abelian Lubin-Tate theory or the conjecture of Deligne-Carayol (\cf \cite{MR1044827}),
and proved by Harris and Taylor \cite{MR1876802}.

The proof of Harris and Taylor was accomplished by global methods;
they related the Lubin-Tate tower to the bad reduction of certain Shimura varieties,
and deduced the non-abelian Lubin-Tate theory from the global theory of Shimura varieties
and automorphic representations.
Although local study of the Shimura varieties at bad primes plays very important role in \cite{MR1876802},
they do not consider the Lubin-Tate tower directly;
their local study is mainly on the Igusa variety, which is in some sense ``vertical'' to the Lubin-Tate tower.

On the other hand, for the action of $\GL_d(F)$ and $D^\times$, there is a purely local study
by Strauch \cite{MR2383890}. Inspired by a pioneering work of Faltings \cite{MR1302321},
he investigated the virtual representation $H_{\mathrm{LT}}$ as a $\GL_d(F)\times D^\times$-module
by using the Lefschetz trace formula, and obtained the result that $H_{\mathrm{LT}}$ realizes
a supercuspidal part of the local Jacquet-Langlands correspondence.
One of nice points of his approach is that one can directly observe that the character relation,
that characterizes the local Jacquet-Langlands correspondence, actually appears in $H_{\mathrm{LT}}$.

Roughly speaking, the work \cite{MR2383890} can be divided into two parts;
the geometric part and the representation-theoretic part.
In the first part (\cite[\S2--3]{MR2383890}), Strauch established the Lefschetz trace formula
for the Lubin-Tate tower and counted the number of fixed points under the group action on the tower.
In the second part, he deduced from the Lefschetz trace formula that the local Jacquet-Langlands correspondence
appears in $H_{\mathrm{LT}}$.
The purpose of this article is to give an alternative approach to the latter, that is,
the representation-theoretic part of \cite{MR2383890}.

The main difference is as follows. In \cite{MR2383890}, Strauch considered the $\pi$-isotypic part $H_{\mathrm{LT}}[\pi]$ of $H_{\mathrm{LT}}$
for a fixed supercuspidal representation $\pi$ of $\GL_d(F)$.
On the other hand, we will begin with an irreducible smooth representation $\rho$ of $D^\times$
and consider the $\rho$-isotypic part $H_{\mathrm{LT}}[\rho]$.
Although the final consequences are more or less similar, our approach has some advantage.
In fact, we can obtain some results on a non-supercuspidal representation of $\GL_d(F)$. This is because the image $\JL(\rho)$ of
$\rho$ under the local Jacquet-Langlands correspondence is not necessarily supercuspidal.
Here is a part of our main results:

\begin{thm}[Theorem \ref{thm:LT-up-to-induction}]\label{thm:Intro}
 As a virtual representation, $H_{\mathrm{LT}}[\rho]-(-1)^{d-1}d\JL(\rho)$ is a linear combination of
 parabolically induced representations. In particular, $H_{\mathrm{LT}}[\rho]$ is not zero.
\end{thm}

For a more precise statement, see Section \ref{sec:LT-rep}.
This theorem implies the following consequence:

\begin{cor}[\cf Corollary \ref{cor:LT-elliptic}]
 If an irreducible discrete series representation $\pi$ appears in $H_{\mathrm{LT}}[\rho]$, then
 $\pi=\JL(\rho)$.
\end{cor}

Together with \cite{non-cusp}, we can also deduce the main result of Strauch \cite[Theorem 4.1.3]{MR2383890}
from Theorem \ref{thm:Intro}:

\begin{cor}[Corollary \ref{cor:LT-cuspidal}]
 For an irreducible supercuspidal representation $\pi$ of $\GL_d(F)$,
 the $\pi$-isotypic part $H^{d-1}_{\mathrm{LT}}[\pi]$ of $H^{d-1}_{\mathrm{LT}}$ is
 isomorphic to $\rho^{\oplus d}$, where $\rho$ is the irreducible smooth representation of $D^\times$ with
 $\pi=\JL(\rho)$.
\end{cor}

The main tool in the argument of \cite[\S 4.1]{MR2383890} is the description of a supercuspidal representation
as the compact induction from a compact-mod-center subgroup, which is a part of the classification result
in \cite{MR1204652}.
On the other hand, in this article we extensively use local harmonic analysis,
such as transfer of orbital integrals. 
The author thinks it possible to extend our method to many other Rapoport-Zink spaces.
For example, in the forthcoming paper \cite{LT-GSp4}, we will investigate the alternating sum of the $\ell$-adic
cohomology of the Rapoport-Zink tower for $\GSp(4)$, by combining the method in this paper with
the Lefschetz trace formula developed in \cite{adicLTF}; in this case, we should replace orbital integrals
by stable orbital integrals.

We sketch the outline of this paper. In Section 2, we briefly recall definitions and notation
concerning the Lubin-Tate tower. In Section 3, we deduce a relation between the distribution character
of $H_{\mathrm{LT}}[\rho]$ and that of $\rho$ from the Lefschetz trace formula.
The author prefers this relation, since it is compatible with the modern characterization of
local endoscopic lifts. In Section 4, we prove that, over the subset of $\GL_d(F)$ consisting of
elliptic regular elements, the distribution character of $H_{\mathrm{LT}}[\rho]$ is
given by a locally constant function $\theta_\rho$, which is strongly related to the (usual) character of $\rho$.
Interestingly, in this section we do not need the character theory due to Harish-Chandra. 
In Section 5, we derive some representation-theoretic conclusions from our character relation.
In this step, in addition to the local Jacquet-Langlands correspondence, we require a rather strong
finiteness property on $H_{\mathrm{LT}}[\rho]$, which is proved by using the Faltings-Fargues isomorphism (\cite{MR1936369}, \cite{MR2441311}).

\bigbreak

\noindent\textbf{Acknowledgment}\quad
The author would like to thank Noriyuki Abe and Matthias Strauch for useful discussions.

\bigbreak

\noindent{\bfseries Notation}
 As in the introduction above, let $F$ be a $p$-adic field and $\mathcal{O}$ its ring of integers.
 We denote the normalized valuation of $F$ by $v_F$.
 For $a\in F$, put $\lvert a\rvert_F=q^{-v_F(a)}$, where $q$ denotes the cardinality of the residue field of $\mathcal{O}$.
 Fix a uniformizer $\varpi$ of $\mathcal{O}$.
 Denote the completion of the maximal unramified extension of $\mathcal{O}$ by $\breve{\mathcal{O}}$ and
 the fraction field of $\breve{\mathcal{O}}$ by $\breve{F}$.
 
 For a reductive algebraic group $G$ over $F$, we sometimes denote $G(F)$ simply by $G$.
 An element $g$ of $G(F)$ is said to be regular if the centralizer $Z(g)$ of $g$ is a maximal torus of $G$.
 In particular, a regular element is semisimple. We denote by $G(F)^\mathrm{reg}$ the set of regular elements in $G(F)$.
 An element $g$ is said to be elliptic if it is contained in an elliptic maximal torus, that is, a maximal torus
 which is anisotropic modulo the center of $G$. A regular element is elliptic if and only if its centralizer is a elliptic maximal torus.
 We denote by $G(F)^\mathrm{ell}$ the set of elliptic regular elements in $G(F)$. Please do not confuse it with the set of elliptic elements.
 
 For a field $k$, we denote its separable closure by $\overline{k}$.

 Let $\ell$ be a prime which is invertible in $\mathcal{O}$. Every representation is considered over $\overline{\Q}_\ell$.

\section{Lubin-Tate tower}\label{sec:Lubin-Tate}
 Let us recall briefly the definition of the Lubin-Tate tower.
 See \cite[\S 2.1]{MR2383890} for more detail.
 Fix an integer $d\ge 1$ and denote by $\mathbb{X}$ a formal $\mathcal{O}$-module over $\overline{\mathbb{F}}_q$
 with $\mathcal{O}$-height $d$ (such $\mathbb{X}$ is unique up to isomorphism).
 For integers $m\ge 0$ and $j$, let $\mathcal{M}^{(j)}_m$ denotes the following functor from the category of complete noetherian local
 $\breve{\mathcal{O}}$-algebras to the category of sets; $\mathcal{M}_m^{(j)}(A)$ consists of isomorphism classes of triples $(X,\rho,\eta)$, 
 where $X$ is a formal $\mathcal{O}$-module over $A$, $\rho$ is an $\mathcal{O}$-quasi-isogeny of $\mathcal{O}$-height $j$ 
 from $\mathbb{X}$ to $X\otimes_A A/\mathfrak{m}_A$ and
 $\eta$ is a Drinfeld $m$-level structure of $X$. This functor is represented by a complete noetherian local ring $R_m^{(j)}$.
 We denote $\Spf R_m^{(j)}$ by $\mathcal{M}^{(j)}_m$ again, and put $\mathcal{M}_m=\coprod_{j\in\Z}\mathcal{M}^{(j)}_m$.

 Put $K_0=\mathrm{GL}_d(\mathcal{O})$ and $K_m=\Ker(\GL_d(\mathcal{O})\to \GL_d(\mathcal{O}/(\varpi^m)))$ for $m\ge 1$.
 We can associate to a compact open subgroup $K$ of $K_0$ the formal scheme $\mathcal{M}_K=\coprod_{j\in\Z}\mathcal{M}^{(j)}_K$,
 so that $\mathcal{M}_{K_m}$ coincides with $\mathcal{M}_m$. Indeed, if $K$ contains $K_m$ for some $m\ge 0$, we define $\mathcal{M}_K$ as the quotient of
 $\mathcal{M}_{K_m}$ by $K/K_m$, where the action of $K$ on $\mathcal{M}_{K_m}$ is given via the Drinfeld $m$-level structures.
 Now we get the projective system of formal schemes $\{\mathcal{M}_K\}_{K\subset K_0}$, which is called the Lubin-Tate tower.

 Let $D$ be the central division algebra over $F$ with invariant $1/d$.
 Then $D^\times$ naturally acts on each $\mathcal{M}_K$, because it is isomorphic to the group of self $\mathcal{O}$-quasi-isogenies of $\mathbb{X}$.
 On the other hand, for $g\in \GL_d(F)$ and a compact open subgroup $K$ of $K_0$ with $g^{-1}Kg\subset K_0$,
 we can define a morphism $\mathcal{M}_K\longrightarrow \mathcal{M}_{g^{-1}Kg}$. Thus we have the action of $\GL_d(F)$ on
 the tower $\{\mathcal{M}_K\}_{K\subset K_0}$ as a pro-object. 

 Denote the rigid generic fiber of $\mathcal{M}_K$ by $M_K$. It is the generic fiber of the adic space $t(\mathcal{M}_K)$ associated to $\mathcal{M}_K$.
 We simply write $H^i_c(M_K)$ for the $\ell$-adic cohomology $H^i_c((M_K)\otimes_{\breve{F}}\overline{\breve{F}},\overline{\Q}_\ell)$,
 and put $H^i_c(M_\infty)=\varinjlim_K H^i_c(M_K)$. It is a smooth representation of $\GL_d(F)\times D^\times$.
 In fact, we may also define the action of the Weil group of $F$ on $H^i_c(M_K)$ and $H^i_c(M_\infty)$, but in this note we do not consider it.

 However this representation is too large to deal with. It is mainly because the formal scheme $\mathcal{M}_K$ is not quasi-compact.
 To make it quasi-compact, we take the quotient by the action of $\varpi^\Z\subset D^\times$ (or $\varpi^\Z\subset \GL_n(F)$, since the action of
 $F^\times\subset \GL_n(F)\times D^\times$, embedded diagonally, is trivial).
 Namely, we shall consider the representation $H^i_c(M_\infty/\varpi^\Z)$ in place of $H^i_c(M_\infty)$.
 This change makes no serious problem;  as in the beginning of the proof of \cite[Theorem 4.1.3]{MR2383890}, 
 we can reduce the latter to the former by twisting.
 Set $H^i_{\mathrm{LT}}:=H^i_c(M_\infty/\varpi^\Z)$ for simplicity.
 
\section{Application of Lefschetz trace formula}

Let us fix an irreducible smooth representation $\rho$ of $D^\times$ which is trivial on $\varpi^\Z\subset D^\times$.
Then $H^i_{\mathrm{LT}}[\rho]:=\Hom_{D^\times}(H^i_{\mathrm{LT}},\rho)^{\mathrm{sm}}$ is an admissible representation of $\GL_n(F)$
(here $(-)^{\mathrm{sm}}$ denotes the subset of smooth vectors).
We would like to study the alternating sum $H_{\mathrm{LT}}[\rho]:=\sum_{i}(-1)^i H^i_{\mathrm{LT}}[\rho]$ by the Lefschetz trace formula.

The following Lefschetz trace formula was established in \cite{MR2383890}. See also \cite{adicLTF}
for more transparent approach to the Lefschetz trace formula for rigid spaces.

\begin{prop}[{{\cite[Theorem 3.3.1]{MR2383890}, \cite[Example 4.21]{adicLTF}}}]\label{prop:LTF}
 Let $(g,b)$ be an element of $\GL_d(F)\times D^\times$ and $K$ a compact open subgroup of $K_0$.
 Assume that $g$ normalizes $K$ and $gK$ consists of elliptic regular elements. Then we have
 \[
 \sum_i(-1)^i\Tr\bigl((g,b);H^i_c(M_K/\varpi^\Z)\bigr)=\#\Fix (g,b).
 \]
 Here the right hand side denotes the number of fixed points under the action of $(g,b)$ on $M_K/\varpi^\Z$.
\end{prop}

Before applying this formula, we will introduce some notation on harmonic analysis.
Let us fix arbitrary Haar measures of $\GL_d(F)$ and $D^\times$.
For $b\in D^\times$, there exists an elliptic element $g_b$ of $\GL_d(F)$ which is stably conjugate to $b$.
Such an element is unique up to conjugacy, and the map $[b]\longmapsto [g_b]$ gives a bijection between the set of elliptic conjugacy classes of $D^\times$
and that of $\GL_d(F)$.
If $b$ is regular, there is an isomorphism $Z(b)\yrightarrow{\cong} Z(g_b)$ of tori which is determined naturally up to conjugacy
(\cf the proof of Lemma \ref{lem:ell-torus-conj}).
We fix Haar measures of all centralizers of $\GL_d(F)$ and $D^\times$,
so that the measures of $Z(b)$ and $Z(g_b)$ are the same under the isomorphism above.

For $\varphi\in C^\infty_c(\GL_d(F))$ and $\gamma\in \GL_d(F)$, we denote its orbital integral 
\[
 \int_{Z(\gamma)\backslash\GL_d(F)}\varphi(g^{-1}\gamma g) dg
\]
by $O_\gamma(f)$. Similarly for $D^\times$.
For $\varphi\in C^\infty_c(\GL_d(F))$ and $\varphi^D\in C^\infty_c(D^\times)$, we say that $\varphi^D$ is a transfer of $\varphi$
if $O_b(f^D)=(-1)^{d-1}O_{g_b}(f)$ for every regular $b\in D^\times$.
Note that the direction of our transfer is the inverse of the usual one.

The following lemma gives the existence of a transfer in some cases.

\begin{lem}
 Let $\varphi$ be an element of $C^\infty_c(\GL_d(F)^{\mathrm{ell}})$. Extending by $0$, we regard $\varphi$ as a function on $\GL_d(F)$.
 Then there exists a transfer $\varphi^D\in C_c^\infty(D^\times)$ of $\varphi$.
\end{lem}

\begin{prf}
 By \cite[Lemma 3.3]{MR700135}, we have only to prove that $\varphi$ is type $E$ (\cf \cite[p.~174, Definition A]{MR700135}). 
 Let $T$ be a maximal torus of $\GL_d(F)$ and $\gamma_0\in T\setminus T^{\mathrm{reg}}$.
 It suffices to observe that $O_\gamma(\varphi)=0$ on some neighborhood of $\gamma_0$.
 We use Shalika's germ expansion (\cf \cite[p.~170]{MR700135}). 
 Let $\{u_j\}$ be a set of representatives for the unipotent conjugacy classes which contain an element commuting with $\gamma_0$,
 such that $u_j$ commutes with $\gamma_0$. Then there exist functions $\Gamma_{\gamma_0,j}^T$ on $T^{\mathrm{reg}}$
 and a neighborhood $N(\varphi)$ of $\gamma_0$ such that
 \[
  O_\gamma(\varphi)=\sum_j O_{\gamma_0 u_j}(\varphi)\Gamma_{\gamma_0,j}^T(\gamma)
 \]
 for all $\gamma\in T^{\mathrm{reg}}\cap N(\varphi)$. Since $\gamma_0 u_j$ is not regular, $\varphi$ vanishes on the conjugacy class of $\gamma_0 u_j$. 
 Thus we have $O_\gamma(\varphi)=0$ for $\gamma\in T^{\mathrm{reg}}\cap N(\varphi)$.
\end{prf}

We endow the discrete group $\varpi^\Z$ with the counting measure.
By using orbital integral, the number of fixed points on $M_K/\varpi^\Z$ is given by the following formula:

\begin{prop}[{{\cite[Theorem 2.6.8]{MR2383890}}}]\label{prop:counting}
 Let $g$ be an element of $\GL_d(F)$ and $b$ an element of $D^\times$. Let $K$ be a compact open subgroup of $K_0$ which is normalized by $g$.
 \begin{enumerate}
  \item Unless $v_F(\det g)-v_F(\Nrd b)$ is divisible by $d$, $\#\Fix(g^{-1},b^{-1})=0$.
  \item Assume one of the following:
	\begin{itemize}
	 \item $v_F(\det g)=v_F(\Nrd b)$ and $b$ is regular,
	 \item or $gK$ consists of regular elements.
	\end{itemize}
	Then $\#\Fix(g^{-1},b^{-1})$ is equal to $d\vol(Z(b)/\varpi^\Z) O_{g_b}(\vol(K)^{-1}\mathbf{1}_{gK})$, where $\mathbf{1}_{gK}$ denotes the
	characteristic function of $gK$.
 \end{enumerate}
\end{prop}

\begin{prf}
 Let $\wp\colon M_K\longrightarrow \mathcal{F}\cong \P^{d-1}$ be the period map. Let us take a fixed point $x\in\mathcal{F}$ under the action of $b^{-1}$,
 and consider fixed points on $\wp^{-1}(x)$ under the action of $(g^{-1},b^{-1})$. In \cite[2.6]{MR2383890}, it is proved that
 such fixed points are in bijection with elements of the set $S_{(g,b)}=\{h\in \GL_d(F)/\varpi^\Z K\mid h^{-1}g_bh\varpi^\Z K=g\varpi^\Z K\}$.

 Assume that $S_{(g,b)}$ is not empty. Then there exist $h\in \GL_d(F)$ and an integer $n$ such that $h^{-1}g_bh\in g\varpi^n K$.
 Then $v_F(\det g_b)=v_F(\det g)+dn$. Since $\det g_b=\Nrd b$, $v_F(\det g)-v_F(\Nrd b)$ is divisible by $d$.
 Thus we have i). Moreover, if $v_F(\det g)=v_F(\Nrd b)$, then $n=0$. In other words, $h$ satisfies $h^{-1}g_bhK=gK$.
 
 Let us prove ii). First assume that $v_F(\det g)=v_F(\Nrd b)$ and $b$ is regular. Then there are $d$ fixed points on $\mathcal{F}$, and we have
 \begin{align*}
 \#\Fix(g^{-1},b^{-1})&=d\#S_{(g,b)}=d\#\{h\in \GL_d(F)/\varpi^\Z K\mid h^{-1}g_bhK=gK\}\\
  &=\frac{d}{\vol(\varpi^\Z K/\varpi^\Z)}\int_{\GL_d(F)/\varpi^\Z}\mathbf{1}_{gK}(h^{-1}g_bh)dh\\
  &=d\cdot \frac{\vol(Z(g_b)/\varpi^\Z)}{\vol(K)}\int_{Z(g_b)\backslash \GL_d(F)}\mathbf{1}_{gK}(h^{-1}g_bh)dh\\
  &=d\vol(Z(g_b)/\varpi^\Z)O_{g_b}\Bigl(\frac{\mathbf{1}_{gK}}{\vol(K)}\Bigr)\\
  &=d\vol(Z(b)/\varpi^\Z)O_{g_b}\Bigl(\frac{\mathbf{1}_{gK}}{\vol(K)}\Bigr).
 \end{align*}
 Next assume that $b$ is not regular and $gK$ consists of regular elements. Then, since $h^{-1}g_bh$ is not regular, it does not belong to $g\varpi^\Z K$.
 Therefore $\#S_{(g,b)}=0$. Similarly the orbital integral $O_{g_b}(\vol(K)^{-1}\mathbf{1}_{gK})$ is $0$, and we get the equality.
\end{prf}

Putting these results together, we obtain the following:

\begin{prop}\label{prop:LT-character-group-element}
 Let $g$ be an element of $\GL_d(F)$ and $K$ a compact open subgroup of $K_0$ such that $K=gKg^{-1}$ and $gK$ consists of elliptic regular elements.
 Put $\varphi_{gK}=\vol(K)^{-1}\mathbf{1}_{gK}$.
 Let $\varphi_{gK}^D\in C^\infty_c(D^\times)$ be a transfer of $\varphi_{gK}$. Then we have
 \[
  \Tr\bigl(\varphi_{gK};H_{\mathrm{LT}}[\rho]\bigr)=(-1)^{d-1}d\Tr(\varphi_{gK}^D;\rho).
 \]
\end{prop}

\begin{prf}
 Put $(D^\times)_g=\{b\in D^\times\mid v_F(\Nrd b)=v_F(\det g)\}$.
 By Proposition \ref{prop:LTF} and Proposition \ref{prop:counting}, we have
 \begin{align*}
  \Tr\bigl(\varphi_{gK};H_{\mathrm{LT}}[\rho]\bigr)&=\Tr\bigl(g;(H_{\mathrm{LT}}[\rho])^K\bigr)
  =\Tr\bigl(g;\Hom_{D^\times}\bigl((H_{\mathrm{LT}})^K,\rho\bigr)\bigr)\\
  &=\Tr\bigl(g;\Hom_{D^\times}\bigl((H^*_c(M_K/\varpi^\Z),\rho\bigr)\bigr)\\
  &=\frac{1}{\vol(D^\times/\varpi^\Z)}\int_{D^\times/\varpi^\Z}\Tr\bigl((g^{-1},b^{-1});H^*_c(M_K/\varpi^\Z)\bigr)\!\Tr(b;\rho)db\\
  &=\frac{1}{\vol(D^\times/\varpi^\Z)}\int_{(D^\times)_g}\#\Fix(g^{-1},b^{-1})\Tr(b;\rho)db\\
  &=d\cdot \frac{\vol\bigl(Z(b)/\varpi^\Z\bigr)}{\vol(D^\times/\varpi^\Z)}\int_{(D^\times)_g}O_{g_b}(\varphi_{gK})\Tr(b;\rho)db\\
  &=(-1)^{d-1}d\cdot \frac{\vol\bigl(Z(b)/\varpi^\Z\bigr)}{\vol(D^\times/\varpi^\Z)}\int_{(D^\times)_g}O_b(\varphi_{gK}^D)\Tr(b;\rho)db\\
  &=(-1)^{d-1}d\cdot \frac{\vol\bigl(Z(b)/\varpi^\Z\bigr)}{\vol(D^\times/\varpi^\Z)}\int_{(D^\times)_g}\int_{Z(b)\backslash D^\times} \varphi_{gK}^D(h^{-1}bh)\Tr(h^{-1}bh;\rho)dhdb\\
  &=\frac{(-1)^{d-1}d}{\vol(D^\times/\varpi^\Z)}\int_{D^\times}\int_{D^\times/\varpi^\Z}\varphi_{gK}^D(h^{-1}bh)\Tr(h^{-1}bh;\rho)dhdb\\
  &=\frac{(-1)^{d-1}d}{\vol(D^\times/\varpi^\Z)}\int_{D^\times/\varpi^\Z}\int_{D^\times}\varphi_{gK}^D(b)\Tr(b;\rho)dbdh\\
  &=\frac{(-1)^{d-1}d}{\vol(D^\times/\varpi^\Z)}\int_{D^\times/\varpi^\Z}\Tr(\varphi_{gK}^D;\rho)dh\\
  &=(-1)^{d-1}d\Tr(\varphi_{gK}^D;\rho),
 \end{align*}
 as desired.
\end{prf}

We want to replace $\varphi_{gK}$ in the proposition above by an arbitrary $\varphi$ which is supported on
$\GL_d(F)^{\mathrm{ell}}$. To do it, we need the following two elementary lemmas:

\begin{lem}\label{lem:elliptic-compact}
 For every elliptic element $g$ of $\GL_d(F)$, we can find a compact open subgroup of $\PGL_d(F)$ which contains the image of $g$ in $\PGL_d(F)$.
\end{lem}

\begin{prf}
 Since $g$ is elliptic, it is contained in a maximal torus which is compact modulo the center of $\GL_d(F)$.
 In particular, its image $\overline{g}$ in $\PGL_d(F)$ topologically generates a compact abelian subgroup $H$.
 Take a compact open subgroup $K$ of $\PGL_d(F)$. Since $H/H\cap K$ is a finite group, we have an integer $n\ge 1$ such that $\overline{g}^n\in K$.
 Therefore $K'=\bigcap_{i\in\Z}\overline{g}^{\,i}K\overline{g}^{\,-i}=\bigcap_{i=0}^{n-1}\overline{g}^{\,i}K\overline{g}^{\,-i}$ is
 a compact open subgroup of $\PGL_d(F)$ which is normalized by $H$.
 Then $HK'=\{hk\mid h\in H, k\in K'\}$ is a compact open subgroup of $\PGL_d(F)$ containing $\overline{g}$.
\end{prf}

\begin{lem}\label{lem:parahoric}
 Let $U$ be the inverse image of a compact open subgroup of $\PGL_d(F)$ under the canonical map $\GL_d(F)\longrightarrow \PGL_d(F)$.
 Then there exists a fundamental system of neighborhoods $\{U_n\}$ of $1\in U$ consisting of compact open normal subgroups of $U$.
\end{lem}

\begin{prf}
 Denote the natural projection $\GL_d(F)\longrightarrow \GL_d(F)/\varpi^\Z$ by $\pr$. Note that $\pr(U)$ is a compact open subgroup of $\GL_d(F)/\varpi^\Z$.

 Let $U'$ be a compact open subgroup of $U$. It suffices to find a compact open subgroup $U''$ of $U'$ which is normal in $U$.
 Since $\pr(U')$ is an open subgroup of a profinite group $\pr(U)$, there exists an open normal subgroup $K$ of $\pr(U)$ contained in $\pr(U')$.
 Put $U''=\pr^{-1}(K)\cap U'$. It is an open subgroup of $U'$, and thus is compact.
 Let us observe that it is normalized by an element $g$ of $U$. For an element $u\in U''$, $\pr(gug^{-1})$ lies in $K=\pr(U'')$.
 In other words, $gug^{-1}\varpi^n\in U''$ for some $n\in\Z$. Since the determinant of every element in $U''$ is a unit of $\mathcal{O}$,
 we may conclude that $n=0$. This completes the proof.
\end{prf}

\begin{thm}\label{thm:LT-character}
 Let $\varphi$ be an element of $C^\infty_c(\GL_d(F)^{\mathrm{ell}})$ and $\varphi^D\in C^\infty_c(D^\times)$ its transfer.
 Then we have
 \[
  \Tr\bigl(\varphi;H_{\mathrm{LT}}[\rho]\bigr)=(-1)^{d-1}d\Tr(\varphi^D;\rho).
 \]
\end{thm}

\begin{prf}
 First note that the right hand side is independent of the choice of a transfer $\varphi^D$, since it depends only on
 the orbital integral of $\varphi^D$ (for example, by Weyl's integral formula).
 Therefore, by Lemma \ref{lem:elliptic-compact}, we may assume that the support of $\varphi$ is contained in some open subgroup $U$
 such as in Lemma \ref{lem:parahoric}.
 Take $\{U_n\}$ as in Lemma \ref{lem:parahoric}. Then $\varphi$ can be written as a linear combination of $\mathbf{1}_{gU_n}$ with
 $g\in U$. Since every maximal compact open subgroup of $\GL_d(F)$ is conjugate to
 $K_0$, there exists $h\in\GL_d(F)$ such that $U'_n=hU_nh^{-1}$ is contained in $K_0$.
 Put $g'=hgh^{-1}$.
 As $g'$ normalizes $U'_n$, by Proposition \ref{prop:LT-character-group-element} we have
 \[
  \Tr\bigl(\mathbf{1}_{g'U'_n};H_{\mathrm{LT}}[\rho]\bigr)=(-1)^{d-1}d\Tr\bigl((\mathbf{1}_{g'U'_n})^D;\rho\bigr).
 \]
 It is clear that the left hand side is equal to $\Tr(\mathbf{1}_{gU_n};H_{\mathrm{LT}}[\rho])$.
 On the other hand, $(\mathbf{1}_{g'U'_n})^D$ also gives a transfer of $\mathbf{1}_{gU_n}$,
 since $\mathbf{1}_{gU_n}$ and $\mathbf{1}_{g'U'_n}$ have the same orbital integrals.
 Therefore the right hand side coincides with $(-1)^{d-1}d\Tr((\mathbf{1}_{gU_n})^D;\rho)$.
 This concludes the proof.
\end{prf}

\section{Character of $H_{\mathrm{LT}}[\rho]$ as a function}
We continue to use the notation in the previous section. 
Let $\Theta_{\mathrm{LT}}$ be the distribution over $\GL_d(F)^{\mathrm{ell}}$ defined by
$\varphi\longmapsto \Tr(\varphi;H_{\mathrm{LT}}[\rho])$. We want to show that it is given by a locally constant function on $\GL_d(F)^{\mathrm{ell}}$
which is strongly related to the character of $\rho$.

Put $\mathbf{ch}=F^{d-1}\times F^\times$. Let $\chi_{\mathrm{GL}_d(F)}\colon \mathrm{GL}_d(F)\longrightarrow \mathbf{ch}$ be
the map that associates to $g\in \GL_d(F)$ the coefficients of the characteristic polynomial of $g$.
Similarly we define $\chi_{D^\times}\colon D^\times\longrightarrow \mathbf{ch}$.
We set $\mathbf{ch}^{\mathrm{ell}}=\chi_{\GL_d(F)}(\GL_d(F)^{\mathrm{ell}})=\chi_{D^\times}((D^\times)^{\mathrm{reg}})$.

\begin{lem}\label{lem:char-proper}
 The map $\chi_{D^\times}$ is proper.
\end{lem}

\begin{prf}
 For an integer $n$, put $\mathbf{ch}_n=\{(a_1,\ldots,a_d)\mid v_F(a_d)=n\}$. Then $\mathbf{ch}$ is the disjoint union of $\mathbf{ch}_n$ as
 a topological space. It suffices to show that $\chi_{D^\times}\colon \chi_{D^\times}^{-1}(\mathbf{ch}_n)\longrightarrow \mathbf{ch}_n$ is proper,
 but it is clear because $\chi_{D^\times}^{-1}(\mathbf{ch}_n)$ is compact and $\mathbf{ch}_n$ is Hausdorff.
\end{prf}

\begin{lem}\label{lem:ell-char-existence}
 There exists a unique locally constant function $\theta_\rho$ on $\GL_d(F)^{\mathrm{ell}}$ such that
 $\theta_\rho(g)=\Tr(b;\rho)$ whenever $g\in \GL_d(F)^{\mathrm{ell}}$ and $b\in D^\times$ are (stably) conjugate.
\end{lem}

\begin{prf}
 We endow $\overline{\Q}_\ell$ with the discrete topology.
 The continuous function $(D^\times)^{\mathrm{reg}}\longrightarrow \overline{\Q}_\ell$; $b\longmapsto \Tr(b;\rho)$ descends to
 a function $\mathbf{ch}^{\mathrm{ell}}\longrightarrow \overline{\Q}_\ell$, which is continuous by Lemma \ref{lem:char-proper}
 (note that the fiber of $\chi_{D^\times}$ over $\mathbf{ch}^{\mathrm{ell}}$ consists of one conjugacy class of $D^\times$).
 Then we may define $\theta_\rho$ as the composite
 $\GL_d(F)^{\mathrm{ell}}\yrightarrow{\chi_{\GL_d(F)}} \mathbf{ch}^{\mathrm{ell}}\longrightarrow \overline{\Q}_\ell$.
\end{prf}

\begin{thm}\label{thm:char-rel}
 We have $\Theta_{\mathrm{LT}}=d\cdot \theta_\rho$ as distributions over $\GL_d(F)^{\mathrm{ell}}$. Namely, for every 
 $\varphi\in C^\infty_c(\GL_d(F)^{\mathrm{ell}})$, we have
 \[
 \Tr\bigl(\varphi;H_{\mathrm{LT}}[\rho]\bigr)=d\int_{\GL_d(F)^{\mathrm{ell}}}\varphi(g)\theta_\rho(g)dg.
 \]
\end{thm}

Let $\varphi^D\in C^\infty_c(D^\times)$ be a transfer of $\varphi$. By Theorem \ref{thm:LT-character}, it suffices to show the equality
\[
 \int_{\GL_d(F)^{\mathrm{ell}}}\varphi(g)\theta_\rho(g)dg=(-1)^{d-1}\Tr(\varphi^D;\rho).
\]
Obviously the right hand side is equal to $(-1)^{d-1}\int_{D^\times}\varphi^D(b)\Tr(b;\rho)db$.
Therefore what we should show is
\[
 \int_{\GL_d(F)^{\mathrm{ell}}}\varphi(g)\theta_\rho(g)dg=(-1)^{d-1}\int_{D^\times}\varphi^D(b)\Tr(b;\rho)db.
\]
The following lemma is used to obtain the equality above:

\begin{lem}\label{lem:ell-torus-conj}
 \begin{enumerate}
  \item There is a natural bijection between the set of conjugacy classes of elliptic maximal tori of $\GL_d(F)$ and 
	the set of conjugacy classes of elliptic maximal tori of $D^\times$.
  \item Let $T$ (resp.\ $T'$) be an elliptic maximal torus of $\GL_d(F)$ (resp.\ $D^\times$) such that $T$ corresponds to $T'$ under the bijection in i).
	Then the rational Weyl group $W_T=N_{\GL_d(F)}(T)/T$ of $T$ and the rational Weyl group $W_{T'}=N_{D^\times}(T')/T'$ of $T'$ 
	are naturally isomorphic.
 \end{enumerate}
\end{lem}

\begin{prf}
 Perhaps this lemma is well-known, but we include its proof for completeness.

 Let us prove i). Let $T$ be an elliptic maximal torus of $\GL_d(F)$. Take $\gamma\in T^{\mathrm{reg}}$ and find $\gamma'\in (D^\times)^{\mathrm{ell}}$
 so that $\gamma'$ is (stably) conjugate to $\gamma$. Let $T'$ be the centralizer of $\gamma'$. Then $T'$ is an elliptic maximal torus of $D^\times$,
 and its conjugacy class depends only on the conjugacy class of $T$.
 Similarly we can attach an elliptic maximal torus of $\GL_d(F)$ to an elliptic maximal torus of $D^\times$. Clearly they are inverse to each other.

 Next consider ii). Let $T$, $T'$, $\gamma$ and $\gamma'$ be as above. 
 Let $W_{T,\overline{F}}^{\Gal}$ be the $\Gal(\overline{F}/F)$-invariant part of the absolute Weyl group 
 $W_{T,\overline{F}}=N_{\GL_d(\overline{F})}(T(\overline{F}))/T(\overline{F})$. Define $W_{T',\overline{F}}^{\Gal}$ similarly.
 Since $H^1(F,T)=H^1(F,T')=0$ by Shapiro's lemma (note that $T$ and $T'$ are the Weil restrictions of $\mathbb{G}_m$ over some finite extensions of $F$),
 the natural inclusions $W_T\hooklongrightarrow W_{T,\overline{F}}^{\Gal}$ and $W_{T'}\hooklongrightarrow W_{T',\overline{F}}^{\Gal}$ are isomorphisms.
 Therefore it suffices to construct an isomorphism from $W_{T,\overline{F}}^{\Gal}$ to $W_{T',\overline{F}}^{\Gal}$.

 Let us fix an isomorphism $\psi\colon \GL_d(\overline{F})\yrightarrow{\cong} (D\otimes_F\overline{F})^\times$ of algebraic group
 over $\overline{F}$ which gives $D^\times$ a structure of an inner form of $\GL_d(F)$.
 Since $\gamma$ and $\gamma'$ are stably conjugate, there exists
 $h\in \GL_d(\overline{F})$ such that $\gamma'=\psi(h\gamma h^{-1})$. Then $\psi_h:=\psi\circ \Ad(h)$ gives an isomorphism $T\yrightarrow{\cong} T'$ over $F$,
 and an isomorphism $W_{T,\overline{F}}\yrightarrow{\cong} W_{T',\overline{F}}$. It suffices to show that
 the second isomorphism induces a bijection from $W_{T,\overline{F}}^{\Gal}$ to $W_{T',\overline{F}}^{\Gal}$.

 Let $n\in N_{\GL_d(\overline{F})}(T(\overline{F}))$ be an element whose image in $W_{T,\overline{F}}$ lies in $W_{T,\overline{F}}^{\Gal}$.
 For $\sigma\in\Gal(\overline{F}/F)$, $\sigma(n)\in nT(\overline{F})$, and thus $\sigma(n\gamma n^{-1})=n\gamma n^{-1}$.
 Namely $n\gamma n^{-1}$ is rational. Hence $\psi_h(n\gamma n^{-1})=\psi_h(n)\gamma'\psi_h(n)^{-1}\in T'(\overline{F})$ is also rational.
 This means that $\psi_h(n)^{-1}\sigma(\psi_h(n))\in T'(\overline{F})$ for every $\sigma\in\Gal(\overline{F}/F)$.
 In other words, the image of $\psi_h(n)$ in $W_{T',\overline{F}}$ is $\Gal(\overline{F}/F)$-invariant.
 Thus $\psi_h$ induces a map from $W_{T,\overline{F}}^{\Gal}$ to $W_{T',\overline{F}}^{\Gal}$.
 Similarly, we can check that $\psi_h^{-1}$ induces a map from $W_{T',\overline{F}}^{\Gal}$
 to $W_{T,\overline{F}}^{\Gal}$. This completes the proof.
\end{prf}

\begin{prf}[of Theorem \ref{thm:char-rel}]
 Let $\{T_i\}_{i\in I}$ (resp.\ $\{T'_i\}_{i\in I}$) be a system of representatives of conjugacy classes of elliptic maximal tori of $\GL_d(F)$
 (resp.\ $D^\times$) such that $T_i$ corresponds to $T_i'$ under the bijection in Lemma \ref{lem:ell-torus-conj} i).
 By Weyl's integral formula, we have
 \begin{align*}
  \int_{\GL_d(F)^{\mathrm{ell}}}\varphi(g)\theta_\rho(g)dg&=\sum_{i\in I}\frac{1}{\#W_{T_i}}\int_{T_i^{\mathrm{reg}}}D_{T_i}(t)\int_{T_i\backslash \GL_d(F)}\varphi(g^{-1}tg)\theta_{\rho}(g^{-1}tg)dgdt\\
  &=\sum_{i\in I}\frac{1}{\#W_{T_i}}\int_{T_i^{\mathrm{reg}}}D_{T_i}(t)\int_{T_i\backslash \GL_d(F)}\varphi(g^{-1}tg)\theta_{\rho}(t)dgdt\\
  &=\sum_{i\in I}\frac{1}{\#W_{T_i}}\int_{T_i^{\mathrm{reg}}}D_{T_i}(t)O_t(\varphi)\theta_{\rho}(t)dt,
 \end{align*}
 where $D_{T_i}(t)=\lvert\det(\Ad(t)-1;\Lie(G)/\Lie(T_i))\rvert_F$ for $t\in T_i^{\mathrm{reg}}$.
 Similarly, we have
 \begin{align*}
  &(-1)^{d-1}\int_{D^\times}\varphi^D(g)\Tr(b;\rho)db\\
  &\qquad=\sum_{i\in I}\frac{(-1)^{d-1}}{\#W_{T'_i}}\int_{T_i'^{\mathrm{reg}}}D_{T'_i}(t')\int_{T_i'\backslash D^\times}\varphi^D(b^{-1}t'b)\Tr(b^{-1}t'b;\rho)dbdt'\\
  &\qquad=\sum_{i\in I}\frac{(-1)^{d-1}}{\#W_{T'_i}}\int_{T_i'^{\mathrm{reg}}}D_{T'_i}(t')\int_{T_i'\backslash D^\times}\varphi^D(b^{-1}t'b)\Tr(t';\rho)dbdt'\\
  &\qquad=\sum_{i\in I}\frac{1}{\#W_{T'_i}}\int_{T_i'^{\mathrm{reg}}}D_{T'_i}(t')\cdot (-1)^{d-1}O_{t'}(\varphi^D)\Tr(t';\rho)dt'.
 \end{align*}
 Let $\psi_i\colon T_i\yrightarrow{\cong}T'_i$ be the isomorphism constructed in the proof in Lemma \ref{lem:ell-torus-conj}.
 It preserves Haar measures on $T_i$ and $T'_i$.
 For $t\in T_i^{\mathrm{reg}}$, put $t'=\psi_i(t)$. Then we have $D_{T_i}(t)=D_{T'_i}(t')$, $O_t(\varphi)=(-1)^{d-1}O_{t'}(\varphi^D)$ and
 $\theta_\rho(t)=\Tr(t';\rho)$. Finally, by Lemma \ref{lem:ell-torus-conj} ii), we have $\#W_{T_i}=\#W_{T'_i}$.
 Therefore we have
 \[
  \int_{\GL_d(F)^{\mathrm{ell}}}\varphi(g)\theta_\rho(g)dg=(-1)^{d-1}\int_{D^\times}\varphi^D(g)\Tr(b;\rho)db,
 \]
 as desired.
\end{prf}

\begin{rem}
 So far, we need neither a finiteness result on $H^i_{\mathrm{LT}}$ nor the deep theory of distribution characters.
\end{rem}

\section{Representations in $H_{\mathrm{LT}}[\rho]$}\label{sec:LT-rep}
In this section, we need the following two deep results:

\begin{thm}\label{thm:LT-fin-length}
 The admissible representation $H^i_{\mathrm{LT}}[\rho]$ has finite length.
\end{thm}

\begin{prf}
 Let $\{M^{\mathrm{Dr}}_n\}_{n\ge 0}$ be the Drinfeld tower. For simplicity we put 
 $H^i_c(M^{\mathrm{Dr}}_n/\varpi^{\Z})=H^i_c((M^{\mathrm{Dr}}_n/\varpi^\Z)\otimes_{\breve{F}}\overline{\breve{F}},\overline{\Q}_\ell)$.
 By the Faltings-Fargues isomorphism (\cite{MR1936369}, \cite{MR2441311}), 
 we have a $\GL_d(F)\times D^\times$-equivariant isomorphism
 $H^i_{\mathrm{LT}}\cong \varinjlim_n H^i_c(M^{\mathrm{Dr}}_n/\varpi^\Z)$. 

 As in \cite[Proposition 4.4.13]{MR2074714}, we can prove that $H^i_c(M^{\mathrm{Dr}}_n/\varpi^\Z)$ is
 a finitely generated $\GL_d(F)$-module.
 Consequently, $(H^i_{\mathrm{LT}})^{K'}$ is a finitely generated $\GL_d(F)$-module for every compact open subgroup $K'\subset D^\times$
 (recall that a submodule of a finitely generated $\GL_d(F)$-module is again finitely generated \cite[Remarque 3.12]{MR771671}).
 Therefore the theorem is reduced to the following general lemma.
\end{prf}

\begin{lem}
 Let $G$ and $H$ be connected reductive algebraic groups over $F$, 
 $\rho$ an irreducible admissible $H$-representation, and $V$ a smooth $G\times H$-representation. Assume the following:
 \begin{itemize}
  \item[(a)] for every compact open subgroup $K'$ of $H$, $V^{K'}$ is a finitely generated $G$-module.
 \end{itemize}
 Then $\Hom_H(\rho,V)$ is a finitely generated smooth $G$-module.
 If moreover 
 \begin{itemize}
  \item[(b)] for every compact open subgroup $K$ of $G$, $V^K$ is a finitely generated $H$-module,
  \item[(c)] and $\rho$ is supercuspidal,
 \end{itemize}
 then $V[\rho]=\Hom_H(V,\rho)^{\mathrm{sm}}$ is an admissible $G$-representation of finite length.
\end{lem}

\begin{prf}
 Since $\rho$ is irreducible, it is generated by one element $x$. Take a compact open subgroup $K'$ of $H$ which stabilizes $x$.
 Then an element $f$ of $\Hom_H(\rho,V)$ is determined by $f(x)$, which lies in $V^{K'}$.
 Therefore we have an injection $\Hom_H(\rho,V)\hooklongrightarrow V^{K'}$ of $G$-modules.
 Since $V^{K'}$ is finitely generated $G$-smooth, so is $\Hom_H(\rho,V)$.

 Next assume also (b) and (c). By (b), $V[\rho]$ is admissible; indeed, $V[\rho]^K=\Hom_H(V^K,\rho)$ for every compact open subgroup $K$ of $G$.
 Let us observe that $\Hom_H(\rho,V)$ is also admissible. Take a compact open subgroup $K$ of $G$ and
 consider the natural $G$-invariant pairing
 \[
  \Hom_H(\rho,V^K)\times \Hom_H(V^K,\rho)\longrightarrow \Hom_H(\rho,\rho)=\overline{\Q}_\ell.
 \]
 By (c), $\rho$ is injective and projective, and thus the pairing above is perfect.
 Therefore we have an injection $\Hom_H(\rho,V)^K=\Hom_H(\rho,V^K)\hooklongrightarrow \Hom_H(V^K,\rho)^*$,
 where $(-)^*$ denotes the (algebraic) dual. Therefore $\Hom_H(\rho,V)^K$ is finite-dimensional,
 and thus $\Hom_H(\rho,V)$ is admissible.
 
 Since $\Hom_H(\rho,V)$ is admissible and finitely generated, it has finite length
 (\cite[Remarque 3.12]{MR771671}).
 The perfect $G$-invariant pairing
 \[
  \Hom_H(\rho,V)\times \Hom_H(V,\rho)\longrightarrow \Hom_H(\rho,\rho)=\overline{\Q}_\ell
 \]
 gives rise to an injection $V[\rho]\hooklongrightarrow \Hom_H(\rho,V)^\vee$, where $(-)^\vee$ denotes
 the contragredient. Now we can conclude that $V[\rho]$ also has finite length.
\end{prf}

\begin{rem}
 The lemma above should be useful in proving the similar result as in Theorem \ref{thm:LT-fin-length} for other Rapoport-Zink spaces,
 once we could establish the analogue of the Faltings-Fargues isomorphism for such Rapoport-Zink spaces.
\end{rem}

\begin{thm}[Local Jacquet-Langlands correspondence]
 Let us denote the set of isomorphism classes of irreducible smooth representations of $D^\times$ by $\mathbf{Irr}(D^\times)$,
 and the set of isomorphism classes of irreducible discrete series representations of $\GL_d(F)$ by $\mathbf{Disc}(\GL_d(F))$.
 There exists a unique map 
 \[
  \JL\colon \mathbf{Irr}(D^\times)\longrightarrow \mathbf{Disc}(\GL_d(F))
 \]
 satisfying the following:
 for an irreducible representation $\rho$ of $D^\times$, the character of $\JL(\rho)$ coincides with $(-1)^{d-1}\theta_\rho$
 over $\GL_d(F)^{\mathrm{ell}}$ (for $\theta_\rho$, see Lemma \ref{lem:ell-char-existence}). 
 The map $\JL$ preserves central characters.
\end{thm}

Let $R(\GL_d(F))$ (resp.\ $R(D^\times)$) be the Grothendieck group over $\overline{\Q}_\ell$ of admissible representations of $\GL_d(F)$ 
(resp.\ $D^\times$) of finite length.
We denote by $R_I(\GL_d(F))$ the $\overline{\Q}_\ell$-subspace of $R(\GL_d(F))$ generated by the images of parabolically induced representations (\cf \cite{MR874042})
and set $\overline{R}(\GL_d(F))=R(\GL_d(F))/R_I(\GL_d(F))$.
For an irreducible representation $\pi$, its image $[\pi]$ in $\overline{R}(\GL_d(F))$ is not equal to $0$ if and only if 
the cuspidal support of $\pi$ is the same as some irreducible discrete series representation
(such $\pi$ is called elliptic, \cf \cite[Lemme 2.1.6]{MR2308851}). Moreover, the local Jacquet-Langlands correspondence
gives an isomorphism $\JL\colon R(D^\times)\yrightarrow{\cong}\overline{R}(\GL_d(F))$ (\cite[Corollaire 2.1.5]{MR2308851}). 

\begin{thm}\label{thm:LT-up-to-induction}
 The image of $H_{\mathrm{LT}}[\rho]$ in $\overline{R}(\GL_d(F))$ is equal to $(-1)^{d-1}d[\JL(\rho)]$.
 In particular, $H_{\mathrm{LT}}[\rho]$ is not zero.
\end{thm}

\begin{prf}
 The image of $H_{\mathrm{LT}}[\rho]$ is written uniquely in the form $\sum_{\tau}a_\tau[\JL(\tau)]$, where 
 $\tau$ runs through isomorphism classes of irreducible smooth representations of $D^\times$ and $a_\tau\in \overline{\Q}_\ell$.
 Since $\varpi^\Z\subset \GL_d(F)$ acts trivially on $H^i_{\mathrm{LT}}[\rho]$, the central character of $\JL(\tau)$ such that $a_\tau\neq 0$ is trivial
 on $\varpi^\Z$. Therefore, if $a_\tau\neq 0$, $\tau$ is trivial on $\varpi^\Z\subset D^\times$.

 By taking characters, we have $d\ch_\rho=\sum_{\tau}(-1)^{d-1}a_\tau\ch_\tau$ over $(D^\times)^{\mathrm{reg}}$
 (\cf Theorem \ref{thm:char-rel}), where $\ch_\tau$ denotes the character of $\tau$.
 By continuity, we have $d\ch_\rho=\sum_{\tau}(-1)^{d-1}a_\tau\ch_\tau$ over $D^\times$.
 By the orthogonality of characters for the compact group $D^\times/\varpi^\Z$, we have
 $a_\tau=(-1)^{d-1}d$ if $\tau\cong \rho$, and $a_\tau=0$ if $\tau\ncong \rho$. This concludes the proof.
\end{prf}

\begin{cor}\label{cor:LT-elliptic}
 Let $\pi$ be an irreducible elliptic representation of $\GL_d(F)$. If the coefficient of $\pi$ in $H_{\mathrm{LT}}[\rho]$ is non-zero,
 then the cuspidal support of $\pi$ is the same as that of $\JL(\rho)$.
\end{cor}

\begin{prf}
 It suffices to note that the image of $\pi$ in $\overline{R}(\GL_d(F))$ is equal to $\pm [\pi']$, 
 where $\pi'$ is the unique discrete series representation whose cuspidal support is the same as that of $\pi$
 (\cf \cite[proof of Lemme 2.1.6]{MR2308851}).
\end{prf}

The following corollary recovers the main result of \cite{MR2383890}:

\begin{cor}\label{cor:LT-cuspidal}
 An irreducible supercuspidal representation $\pi$ of $\GL_d(F)$ which is trivial on $\varpi^\Z\subset \GL_d(F)$
 appears in $H^{d-1}_{\mathrm{LT}}[\rho]$ if and only if $\pi=\JL(\rho)$.
 Moreover, if $\JL(\rho)$ is supercuspidal, then we have $H^{d-1}_{\mathrm{LT}}[\JL(\rho)]=\rho^{\oplus d}$.
\end{cor}

\begin{prf}
 By \cite{non-cusp}, a supercuspidal representation appears only in $H^{d-1}_{\mathrm{LT}}$. Therefore the multiplicity of $\pi$
 in $H^{d-1}_{\mathrm{LT}}[\rho]$ coincides with the coefficient of $\pi$ in $H_{\mathrm{LT}}[\rho]$ multiplied by $(-1)^{d-1}$.
 By Corollary \ref{cor:LT-elliptic}, if it is non-zero then the cuspidal support of $\pi$ coincides with that of $\JL(\rho)$.
 Namely, $\pi=\JL(\rho)$. 

 Moreover, if $\pi=\JL(\rho)$, the multiplicity of $\pi$ in $H^{d-1}_{\mathrm{LT}}[\rho]$ is $d$ by Theorem \ref{thm:LT-up-to-induction}. Therefore the multiplicity of $\rho$ in
 $H^{d-1}_{\mathrm{LT}}[\pi]$ is also $d$. On the other hand, any irreducible smooth representation $\rho'$ of $D^\times$
 which is not isomorphic to $\rho$ does not appear in $H^{d-1}_{\mathrm{LT}}[\pi]$, as $\pi\neq \JL(\rho')$.
 This concludes $H^{d-1}_{\mathrm{LT}}[\JL(\rho)]=\rho^{\oplus d}$, as desired.
\end{prf}

\def\cprime{$'$} \def\cprime{$'$}
\providecommand{\bysame}{\leavevmode\hbox to3em{\hrulefill}\thinspace}
\providecommand{\MR}{\relax\ifhmode\unskip\space\fi MR }
\providecommand{\MRhref}[2]{%
  \href{http://www.ams.org/mathscinet-getitem?mr=#1}{#2}
}
\providecommand{\href}[2]{#2}


\begin{thebibliography}{Mie10b}

\bibitem[Ber84]{MR771671}
J.~N. Bernstein, \emph{Le ``centre'' de {B}ernstein}, Representations of
  reductive groups over a local field, Travaux en Cours, Hermann, Paris, 1984,
  Edited by P. Deligne, pp.~1--32.

\bibitem[BK93]{MR1204652}
Colin~J. Bushnell and Philip~C. Kutzko, \emph{The admissible dual of {${\rm
  GL}(N)$} via compact open subgroups}, Annals of Mathematics Studies, vol.
  129, Princeton University Press, Princeton, NJ, 1993.

\bibitem[Car90]{MR1044827}
H.~Carayol, \emph{Nonabelian {L}ubin-{T}ate theory}, Automorphic forms,
  {S}himura varieties, and {$L$}-functions, {V}ol.\ {II} ({A}nn {A}rbor, {MI},
  1988), Perspect. Math., vol.~11, Academic Press, Boston, MA, 1990,
  pp.~15--39.

\bibitem[Dat07]{MR2308851}
J.-F. Dat, \emph{Th\'eorie de {L}ubin-{T}ate non-ab\'elienne et
  repr\'esentations elliptiques}, Invent. Math. \textbf{169} (2007), no.~1,
  75--152.

\bibitem[Fal94]{MR1302321}
G.~Faltings, \emph{The trace formula and {D}rinfel\cprime d's upper halfplane},
  Duke Math. J. \textbf{76} (1994), no.~2, 467--481.

\bibitem[Fal02]{MR1936369}
\bysame, \emph{A relation between two moduli spaces studied by {V}. {G}.
  {D}rinfeld}, Algebraic number theory and algebraic geometry, Contemp. Math.,
  vol. 300, Amer. Math. Soc., Providence, RI, 2002, pp.~115--129.

\bibitem[Far04]{MR2074714}
L.~Fargues, \emph{Cohomologie des espaces de modules de groupes
  {$p$}-divisibles et correspondances de {L}anglands locales}, Ast\'erisque
  (2004), no.~291, 1--199, Vari{\'e}t{\'e}s de Shimura, espaces de
  Rapoport-Zink et correspondances de Langlands locales.

\bibitem[FGL08]{MR2441311}
L.~Fargues, A.~Genestier, and V.~Lafforgue, \emph{L'isomorphisme entre les
  tours de {L}ubin-{T}ate et de {D}rinfeld}, Progress in Mathematics, vol. 262,
  Birkh\"auser Verlag, Basel, 2008.

\bibitem[HT01]{MR1876802}
M.~Harris and R.~Taylor, \emph{The geometry and cohomology of some simple
  {S}himura varieties}, Annals of Mathematics Studies, vol. 151, Princeton
  University Press, Princeton, NJ, 2001, With an appendix by Vladimir G.
  Berkovich.

\bibitem[Kaz86]{MR874042}
David Kazhdan, \emph{Cuspidal geometry of {$p$}-adic groups}, J. Analyse Math.
  \textbf{47} (1986), 1--36.

\bibitem[Mie]{LT-GSp4}
Y.~Mieda, \emph{Lefschetz trace formula and {$\ell$}-adic cohomology of
  {R}apoport-{Z}ink tower for {${\rm GSp}(4)$}}, in preparation.

\bibitem[Mie10a]{adicLTF}
\bysame, \emph{Lefschetz trace formula for open adic spaces}, preprint,
  arXiv:1011.1720, 2010.

\bibitem[Mie10b]{non-cusp}
\bysame, \emph{Non-cuspidality outside the middle degree of $\ell$-adic
  cohomology of the {L}ubin-{T}ate tower}, Adv. Math. \textbf{225} (2010),
  no.~4, 2287--2297.

\bibitem[Rog83]{MR700135}
J.~D. Rogawski, \emph{Representations of {${\rm GL}(n)$} and division algebras
  over a {$p$}-adic field}, Duke Math. J. \textbf{50} (1983), no.~1, 161--196.

\bibitem[Str08]{MR2383890}
M.~Strauch, \emph{Deformation spaces of one-dimensional formal modules and
  their cohomology}, Adv. Math. \textbf{217} (2008), no.~3, 889--951.

\end{thebibliography}
\end{document}